\documentclass[a4paper,11pt]{article}

\usepackage{graphicx} 
\usepackage{amsmath}
\usepackage{amssymb}
\usepackage{amsfonts}
\usepackage{hyperref} 
\usepackage{amsthm}
\usepackage{mathtools}

\usepackage{tikz}
\usetikzlibrary{graphs,graphs.standard}

\usepackage[utf8]{inputenc}
\usepackage{graphicx}

\usepackage{float}
\usepackage{caption}
\usepackage{diagbox} 

\usepackage[verbose]{placeins}

\usepackage[numbers]{natbib}
\usepackage{cite}

\usepackage{setspace}

\theoremstyle{plain} 
\newtheorem{thm}{Theorem}[section]

\theoremstyle{definition} 

\theoremstyle{remark} 

\title{Enumeration of $n$-plexes}
\author{Arjun Maniyar}
\date{}

\begin{document}

\maketitle

\begin{abstract}
    Palmer provides a method of enumerating \textit{n-plexes} in \citep{PALMER}, however it has some typographical errors in the formula for the cycle index $Z(S_p^{(r)})$ and the values of $s_p^n$, the number of $n$-plexes on $p$ points. This article is intended to provide the correct formulas. 
\end{abstract}

\section{Introduction}

An \textit{$n$-plex} of order $p$ is an $n$-dimensional simplicial complex with $p$ points where every maximal simplex has dimension $n$ or $0$ (a maximal simplex is a simplex that is not a face of any larger simplex). For example, every graph is a $1$-plex. This article is concerned with counting these structures. Note that this problem is the same as counting $(n+1)$-uniform hypergraphs on $p$ vertices. All the ideas and notation used in this article are from \citep{PALMER}, the sole purpose of this article is to provide the correct formulas for $Z(S_p^{(r)})$ and correct values of $s_p^n$.

Let $s_p^n(x)$ be the counting polynomial for $n$-plexes of order $p$, so by definition we have \[s_p^n(x)=\sum s_{p,k}^n x^k,\] where $s_{p,k}^n$ is the number of $n$-plexes with $k$ simplexes of dimension $n$.

In order to count these structures we first need some definitions. Let $A$ be a permutation group acting on a set $X = \{1,2,\dots , p\}$. First recall that we can write any permutation
as a unique product of disjoint cycles.
Now, for a permutation $\alpha$ in $A$, we let \textit{$j_k(\alpha)$} be the number of cycles with exactly $k$ elements. Letting $y_1, y_2, \dots, y_p$ be variables, we define the \textit{cycle index} of $A$ by

    $$ Z(A) = \frac{1}{|A|} \sum_{\alpha \in A} \prod_{k=1}^{p} y_k^{j_k(\alpha)} 
 = \frac{1}{|A|} \sum_{\alpha \in A} y_1^{j_1(\alpha)} \cdots y_p^{j_p(\alpha)}.$$
To display the variables we may also write $Z(A) = Z(A; y_1, \dots, y_p)$. We may also substitute a function $f(x)$ into $Z(A)$, using this notation
\[ Z(A, f(x)) = Z(A; f(x), f(x^2), \dots, f(x^p)). \]

The action of $A$ on the set $X$ induces a corresponding permutation group $A^{(r)}$, called the \textit{$r$-group} of $A$, which acts on $X^{(r)}$ the set of all $r$-element subsets of $X$. That is, each permutation $\alpha \in A$ induces a permutation $\alpha' \in A^{(r)}$ such that for each $\{i_1, \dots, i_r\} \in X^{(r)}$,
\begin{equation} \label{induced perm def}
    \alpha' \{i_1, \dots, i_r\} = \{\alpha i_1, \alpha i_2, \dots, \alpha i_r\}.
\end{equation}
    
We assume that $r\leq p$, so that $X^{(r)}$ is non empty. Note that $A^{(r)}$ and $A^{(p-r)}$ are identical permutation groups. This is because every element in $X^{(r)}$ is in bijective correspondence to an element in $X^{(p-r)}$ by taking complements of the subsets in $X$, and so $A^{(r)}$ and $A^{(p-r)}$ are the same permutation groups upto renaming the elements. Therefore, we have
\begin{equation}\label{eq:cycle index same for r and p-r}
    Z(A^{(r)}) = Z(A^{(p-r)}).
\end{equation}

Frank Harary showed in \citep{harary1955} that P\'olya's enumeration theorem \citep{polya1937} can be used to express the counting polynomial $s_p^1(x)$ for graphs in terms of the cycle index of the pair group $S_p^{(2)}$ of the symmetric group $S_p$. Similarly the following result can also be proved (see \citep[pp.~25, 32]{harary1967} for a proof).

\begin{thm}\label{thm:n-plexes}
    Let $p\geq n+1$ then the counting polynomial $s_p^n$ for $n$-plexes of order $p$ is given by
    $$s_p^n(x) = Z(S_p^{(n+1)}, 1+x),$$
    where $S_p^{(n+1)}$ is the $(n+1)$-group of the symmetric group $S_p$.
\end{thm}
The total number of $n$-plexes $s_p^n$ can be calculated by setting $x=1$ in this formula. To apply this result, the cycle index $Z(S_p^{(n+1)})$ is required. 
The next section deals with deriving a method to compute this cycle index.

\section{Computing the cycle index of $S_p^{(r)}$ }

Let us first find the cycle index of $S_p$. Observe that the cycle index only depends on the length of the cycles in a permutation, not on what's inside the cycles. So, we can group up terms with the same cycle structure. Let $(j) = (j_1, \dots , j_p)$ represent an integer partition of $p$, where $j_k$ is the number of parts equal to $k$. Note that each permutation $\alpha$ in $S_p$ corresponds to an integer partition $(j)$, where $j_k = j_k(\alpha)$ for all $k$. And, since there are $h(j) = p!/(\prod_{k} k^{j_k} j_k!)$ permutations in $S_p$ that correspond to the partition $(j)$, we get the following formula for $Z(S_p)$.
$$Z(S_p) = \frac{1}{p!} \sum_{(j)}  h(j) \prod_{k} y_k^{j_k} = \frac{1}{p!} \sum_{(j)} \frac{p!}{\prod_{k} k^{j_k} j_k!}  y_1^{j_1} \cdots y_p^{j_p},$$
where the sum is over all partitions $(j)$ of $p$.

An explicit formula for $Z(S_p^{(2)})$ for counting graphs is also known and can be found in many books and articles, for instance see \citep[\S ~ 6]{harary1967}. But finding an explicit formula for $Z(S_p^{(r)})$ for $r \geq 3$ gets more difficult and cumbersome, see \citep{PALMER} for an approach to find an explicit formula for $Z(S_p^{(3)})$. However, Palmer provides an easier way in \citep{PALMER} to find the cycle index $Z(S_p^{(r)})$ without finding an explicit formula. This method is outlined below.

Let $(j)$ be an integer partition of $p$ and suppose $\alpha$ is a permutation in $S_p$ whose disjoint cycle decomposition corresponds to $(j)$, that is, $j_k = j_k(\alpha)$ for all $k$. We wish to find the cycle structure of the induced  permutation $\alpha'$ in $S_p^{(r)}$, that is, we want to find $j_k(\alpha')$ for each $k$. First we use the following formula:
\begin{equation}\label{j1}
    j_1(\alpha') = \sum_{(i)} \prod_k \binom{j_k}{i_k} \, ,
\end{equation}
where the sum is over all partitions $(i) = (i_1, i_2, \dots)$ of $r$. To see why this is true consider the case $r=4$. 
Suppose $p=8$ and $\alpha = (1)(2\ 3)(4\ 5)(6\ 7\ 8)$. Now by the definition of $\alpha'$(see equation (\ref{induced perm def})), a 4-subset is fixed by $\alpha'$ if and only if it is a union of disjoint cycles of $\alpha$ whose total length is 4.
In our example, $\alpha'$ only fixes the following two 4-subsets: $\{1,6,7,8\}$ (union of one 3-cycle and one 1-cycle) and $\{2,3,4,5\}$ (union of two 2-cycles). So, in this case $j_1(\alpha') = 2$. 
Similarly for any $\alpha'$ in $S_p^{(4)}$, the number of 1-cycles $j_1(\alpha')$ is given by the following formula
$$\binom{j_1}{4} + \binom{j_1}{2}\binom{j_2}{1} + \binom{j_2}{2} + \binom{j_1}{1}\binom{j_3}{1} + \binom{j_4}{1} \, .$$
Note that if $\alpha_1$ and $\alpha_2$ both have the same cycle structure then $j_1(\alpha_1') = j_1(\alpha_2')$. So, computing $j_1(\alpha')$ for each partition $(j)$ using equation (\ref{j1}), gives us $j_1(\alpha')$ for all permutations $\alpha'$ in $S_p^{(r)}$. 

Now let $(j)$ be a partition of $p$ and $\beta$ be a permutation in $S_p$ associated to this partition. From equation (\ref{j1}), we already have $j_1((\beta^m)')$ for each $m$. Also observe that $(\beta^m)' = (\beta')^m$, this follows directly from the definition of induced permutation (see equation(\ref{induced perm def})). 
Now notice that a $d$-cycle in $\beta'$ becomes a 1-cycle in $(\beta')^m$ if and only if $d|m$, since raising a permutation to $m$-th power fixes exactly those elements whose cycle length divides $m$.
Moreover, if $d|m$, then each $d$-cycle in $\beta'$ becomes $d$ 1-cycles in $(\beta')^m$. So we get the following formula \footnote{In \citep{PALMER}, Palmer includes a factor of $1/m$ in front of the sum in formula (\ref{main formula in the method}), but this appears to be a typo.}
\begin{equation}\label{main formula in the method}
    j_1((\beta')^m) =  \sum_{d|m} d j_d(\beta').
\end{equation}
We can use this formula successively with $m=2,3,4,\dots$ to calculate $j_m(\beta')$. Hence we can determine the contribution of $\beta'$ in $Z(S_p^{(r)})$. This process can be repeated for each partition $(j)$ (we apply this method to only one permutation of a selected partition $(j)$ since all other permutations corresponding to $(j)$ have the same contribution). We finally obtain $Z(S_p^{(r)})$ by adding up all the contributions. We can also use the following formula from \citep{PALMER} to speed up the calculations: if $\beta'$ contributes $\Pi y_k^{j_k}$ to the cycle index, then the contribution of any power $(\beta')^m$ is 
$$\prod_k y_{k/(m,k)}^{(m,k)j_k}.$$

See Section \ref{sec:formulas} for the formulas of the cycle index which were calculated using this method.

\section{Number of $n$-plexes}

In this section we use formulas from Section \ref{sec:formulas} to calculate $s_p^n$, number of $n$-plexes on $p$ points. From Theorem \ref{thm:n-plexes} we already know that 
$$s_p^n(x) = Z(S_p^{(n+1)}, 1+x).$$
To calculate $s_p^n$ we need to substitute $x=1$ in the equation. This is the same as directly substituting $y_i = 2$ for each variable $y_i$ in the cycle index  $Z(S_p^{(n+1)})$. Using the formulas from Section \ref{sec:formulas} then gives us the following table which includes $s_p^n$ for $p\leq 9$ and $n\leq3$.

\begin{table}[h!]
\centering
\caption{The number of \textit{n}-plexes of order $p\leq9$.}
\begin{tabular}{| l | c c c |}
\hline
\diagbox{$p$}{$n$} & 1 & 2 & 3 \\
\hline
1 & 1 & 1 & 1 \\
2 & 2 & 1 & 1 \\
3 & 4 & 2 & 1 \\
4 & 11 & 5 & 2 \\
5 & 34 & 34 & 6\\
6 & 156 & 2136 & 156\\
7 & 1044 & 7013320 & 7013320\\
8 & 12346 & 1788782616656 & 29281354514767168\\
9 & 274668 & 53304527811667897248 & 234431745534048922731115555415680\\
\hline
\end{tabular}
\end{table}

\section{Formulas for $Z(S_p^{(r)})$} \label{sec:formulas}

In this section we list all the formulas for $Z(S_p^{(r)})$ for $p\leq9$. Since exact formulas of $Z(S_p)$ and $Z(S_p^{(2)})$ are already known, so we only provide formulas for the case $r = 3,4$. Using equation (\ref{eq:cycle index same for r and p-r}), the following 6 formulas completes the list for $Z(S_p^{(r)})$ for $p\leq9$. Although we expect $Z(S_p^{(r)})$ to have same number of terms as number of integer partitions of $p$, the expressions below may contain fewer terms since identical terms have been merged. For example, the unmerged cycle index of $S_6^{(3)}$ is:

\begin{align*}
    Z(S_6^{(3)}) = \frac{1}{6!} & (a_{1}^{20} + 15 a_{1}^{8} a_{2}^{6} + 45 a_{1}^{4} a_{2}^{8} + 40 a_{1}^{2} a_{3}^{6} + 40 a_{1}^{2} a_{3}^{6} + 120 a_{1}^{2} a_{3}^{2} a_{6}^{2} \\& + 15 a_{2}^{10} 
    + 90 a_{2}^{2} a_{4}^{4} + 90 a_{2}^{2} a_{4}^{4} + 120 a_{2} a_{6}^{3} + 144 a_{5}^{4}).
\end{align*}
The merged cycle indexes are listed below.

\begin{align*}
    Z(S_6^{(3)}) = \frac{1}{6!} & (a_{1}^{20} + 15 a_{1}^{8} a_{2}^{6} + 45 a_{1}^{4} a_{2}^{8} + 80 a_{1}^{2} a_{3}^{6} + 120 a_{1}^{2} a_{3}^{2} a_{6}^{2} + 15 a_{2}^{10} \\&
    + 180 a_{2}^{2} a_{4}^{4}  + 120 a_{2} a_{6}^{3} + 144 a_{5}^{4}).
\end{align*}

\begin{align*}
    Z(S_7^{(3)}) = & \frac{1}{7!}  (a_{1}^{35} + 21 a_{1}^{15} a_{2}^{10} + 105 a_{1}^{7} a_{2}^{14} + 70 a_{1}^{5} a_{3}^{10} + 105 a_{1}^{3} a_{2}^{16} \\&
    + 420 a_{1}^{3} a_{2} a_{3}^{4} a_{6}^{3} + 280 a_{1}^{2} a_{3}^{11} + 420 a_{1}  a_{4} a_{6} a_{12}^{2}+ 840 a_{1} a_{2}^{3} a_{4}^{7} \\&
    + 210 a_{1} a_{2}^{2}a_{3}^{2} a_{6}^{4} + 504  a_{5}^{3} a_{10}^{2} + 840 a_{2} a_{3} a_{6}^{5} + 504 a_{5}^{7} + 720 a_{7}^{5}).
\end{align*}

\begin{align*}
    Z(S_8^{(3)}) = &  \frac{1}{8!}  (a_{1}^{56} + 28 a_{1}^{26} a_{2}^{15} + 210 a_{1}^{12} a_{2}^{22} + 112 a_{1}^{11} a_{3}^{15} + 420 a_{1}^{6} a_{2}^{25} \\&
    + 1120 a_{1}^{5} a_{2}^{3} a_{3}^{7} a_{6}^{4} + 420 a_{1}^{4} a_{2}^{4} a_{4}^{11} + 1680 a_{1}^{3} a_{2}^{4} a_{3}^{3} a_{6}^{6} + 2520 a_{1}^{2} a_{2}^{5} a_{4}^{11} \\&
    + 1120 a_{1}^{2} a_{3}^{18} + 1120 a_{1}^{2} a_{3}^{8} a_{6}^{5} + 4032 a_{1}  a_{5}^{5} a_{10}^{3}+ 3360 a_{1} a_{2} a_{3} a_{4}^{2} a_{6} a_{12}^{3} \\&
    + 2688 a_{1}  a_{5}^{2} a_{15}^{3} + 1344 a_{1} a_{5}^{11} + 105 a_{2}^{28} + 1260 a_{2}^{6} a_{4}^{11} + 3360 a_{2} a_{3}^{2} a_{6}^{8} \\&
    + 3360 a_{2} a_{6}^{9} + 1260 a_{4}^{14} + 5760 a_{7}^{8} + 5040 a_{8}^{7}).
\end{align*}

\begin{align*}
    Z(S_9^{(3)})  = &  \frac{1}{9!} (a_{1}^{84} + 36 a_{1}^{42} a_{2}^{21} + 168 a_{1}^{21} a_{3}^{21} + 378 a_{1}^{20} a_{2}^{32} + 1260 a_{1}^{10} a_{2}^{37}\\&
    + 756 a_{1}^{10} a_{2}^{5} a_{4}^{16} + 2520 a_{1}^{9} a_{2}^{6} a_{3}^{11} a_{6}^{5} + 7560 a_{1}^{5} a_{2}^{8} a_{3}^{5} a_{6}^{8} + 945 a_{1}^{4} a_{2}^{40} \\&
    + 7560 a_{1}^{4} a_{2}^{8} a_{4}^{16} + 3024 a_{1}^{4} a_{5}^{16} + 5600 a_{1}^{3} a_{3}^{27} + 10080 a_{1}^{3} a_{3}^{13} a_{6}^{7} \\&
    + 18144 a_{1}^{2}  a_{2} a_{5}^{8} a_{10}^{4}+ 11340 a_{1}^{2} a_{2}^{9} a_{4}^{16} + 15120 a_{1}  a_{2}^{2} a_{3}^{3} a_{4}^{4} a_{6} a_{12}^{4} \\&
    + 15120 a_{1}  a_{2}^{2} a_{3} a_{4}^{4} a_{6}^{2} a_{12}^{4}+ 24192 a_{1}  a_{3} a_{5}^{4} a_{15}^{4} + 2520 a_{1} a_{2}^{10} a_{3}^{3} a_{6}^{9} \\&
    + 30240 a_{1} a_{2} a_{3}^{3} a_{6}^{12} + 30240 a_{1} a_{2} a_{3} a_{6}^{13} + 9072  a_{2}^{2} a_{5}^{4} a_{10}^{6}\\&
    + 18144  a_{4} a_{5}^{2} a_{10} a_{20}^{3} + 25920  a_{7}^{6} a_{14}^{3} + 11340 a_{2}^{2} a_{4}^{20} + 40320 a_{3} a_{9}^{9} \\&
    + 45360 a_{4} a_{8}^{10} + 25920 a_{7}^{12}).
\end{align*}

\begin{align*}
    Z(S_8^{(4)}) = &  \frac{1}{8!} (a_{1}^{70} + 28 a_{1}^{30} a_{2}^{20} + 210 a_{1}^{14} a_{2}^{28} + 112 a_{1}^{10} a_{3}^{20} + 525 a_{1}^{6} a_{2}^{32} \\&
    + 1120 a_{1}^{6} a_{2}^{2} a_{3}^{8} a_{6}^{6} + 1120 a_{1}^{4} a_{3}^{22} + 3360 a_{1}^{2}  a_{4}^{2} a_{6}^{2} + 4200 a_{1}^{2} a_{2}^{6} a_{4}^{14} \\&
    + 1680 a_{1}^{2} a_{2}^{4} a_{3}^{4} a_{6}^{8} + 1260 a_{1}^{2} a_{2}^{2} a_{4}^{16} + 4032  a_{5}^{6} a_{10}^{4} + 2688  a_{5}^{2} a_{15}^{4} \\&
    + 1120 a_{2}^{2} a_{3}^{10} a_{6}^{6} + 6720 a_{2}^{2} a_{3}^{2} a_{6}^{10} + 5040 a_{2} a_{4} a_{8}^{8} + 1344 a_{5}^{14} + 5760 a_{7}^{10}).
\end{align*}

\begin{align*}
    Z(S_9^{(4)})  = &  \frac{1}{9!} (a_{1}^{126} + 36 a_{1}^{56} a_{2}^{35} + 378 a_{1}^{26} a_{2}^{50} + 168 a_{1}^{21} a_{3}^{35} + 1260 a_{1}^{12} a_{2}^{57}\\&
    + 2520 a_{1}^{11} a_{2}^{5} a_{3}^{15} a_{6}^{10} + 945 a_{1}^{6} a_{2}^{60} + 756 a_{1}^{6} a_{2}^{10} a_{4}^{25} + 3360 a_{1}^{6} a_{3}^{40} \\&
    + 7560 a_{1}^{5} a_{2}^{8} a_{3}^{7} a_{6}^{14} + 7560 a_{1}^{4} a_{2}^{11} a_{4}^{25} + 15120 a_{1}^{3}  a_{2} a_{3} a_{4}^{4} a_{6}^{3} a_{12}^{7} \\&
    + 2520 a_{1}^{3} a_{2}^{9} a_{3}^{3} a_{6}^{16} + 11340 a_{1}^{2} a_{2}^{12} a_{4}^{25} + 10080 a_{1}^{2} a_{2}^{2} a_{3}^{18} a_{6}^{11}\\&
    + 11340 a_{1}^{2} a_{2}^{2} a_{4}^{30} + 9072 a_{1} a_{5}^{5} a_{10}^{10} + 18144 a_{1}  a_{5}^{11} a_{10}^{7}\\&
    + 18144 a_{1} a_{5}  a_{10}^{2} a_{20}^{5} + 15120 a_{1}  a_{2}^{2} a_{3} a_{4}^{4} a_{6}^{3} a_{12}^{7} + 24192 a_{1}  a_{5}^{4} a_{15}^{7} \\& 
    + 3024 a_{1} a_{5}^{25} + 25920  a_{7}^{8} a_{14}^{5} + 10080 a_{2}^{3} a_{3}^{4} a_{6}^{18} + 30240 a_{2}^{3} a_{3}^{2} a_{6}^{19}\\& 
    + 45360 a_{2} a_{4} a_{8}^{15} + 2240 a_{3}^{42} + 20160 a_{3}^{4} a_{6}^{19}
    + 25920 a_{7}^{18} + 40320 a_{9}^{14}).
\end{align*}

\section*{Acknowledgements} This work was part of my summer project supported by Undergraduate Research Support Scheme (URSS) that I carried out at Warwick Mathematics Institute under the supervision of Helena Verrill. Errors in these formulas were pointed out to me by my supervisor. I am very thankful to my supervisor, Helena Verrill, for their continuous guidance, support, and feedback throughout the course of this work.

\addcontentsline{toc}{section}{References}
\bibliography{biblio}
\end{document}